\RequirePackage{fix-cm}
\documentclass[smallextended]{svjour3}       
\smartqed  
\usepackage{graphicx}
\usepackage{amsmath}
\usepackage{amssymb}
\usepackage{color}
\usepackage{algorithm}
\usepackage{subfig}

\begin{document}

\title{Structure-Preserving Model-Reduction of Dissipative Hamiltonian Systems
}


\author{Babak Maboudi Afkham         \and
        Jan S. Hesthaven 
}


\institute{Babak Maboudi Afkham \at
              EPFL SB MATH MCSS, MA C1 645 (B\^atiment MA), Station 8, CH-1015 Lausanne, Switzerland \\
              Tel.: +41-21-6935905 \\
              \email{babak.maboudi@epfl.ch}           
           \and
           Jan S. Hesthaven \at
           EPFL SB MATH MCSS, MA C2 652 (B\^atiment MA), Station 8, CH-1015 Lausanne, Switzerland
}

\date{Received: date / Accepted: date}

\maketitle

\begin{abstract}
Reduced basis methods are popular for approximately solving large and complex systems of differential equations. However, conventional reduced basis methods do not generally preserve conservation laws and symmetries of the full order model. Here, we present an approach for reduced model construction, that preserves the symplectic symmetry of dissipative Hamiltonian systems. The method constructs a closed reduced Hamiltonian system by coupling the full model with a canonical heat bath. This allows the reduced system to be integrated with a symplectic integrator, resulting in a correct dissipation of energy, preservation of the total energy and, ultimately, in the stability of the solution. Accuracy and stability of the method are illustrated through the numerical simulation of the dissipative wave equation and a port-Hamiltonian model of an electric circuit.

\end{abstract}

\section{Introduction} \label{sec:1}

The need for increased accuracy has led to more complex models and the use of large systems of partial differential equations in engineering and science. As a consequence, direct numerical methods for solving PDEs have become computationally demanding and, at times impractical.

During the past decade, reduced basis methods have emerged as a powerful approach to reduce the cost of evaluating large systems of partial differential equations \cite{Ito:1998up,Ito:1998ch,Ito:2001ev}. These methods construct a low-dimensional linear subspace, the reduced space, that approximately represents the solution to the system of differential equations. The projection of the original system onto the reduced space then allows the exploration of the solution with a significantly reduced computational complexity \cite{hesthaven2015certified,quarteroni2015reduced}.

Hamiltonian systems are an important class of systems that appear in engineering and science. In these systems, preserving system energy is essential to obtain a correct numerical solution. Therefore, the development of model reduction techniques that preserve the symplectic symmetry is crucial. However, the classical reduced basis methods do not generally preserve the conservation laws and intrinsic symmetries of Hamiltonian systems \cite{Amsallem:2014ef,prajna2003pod}. This often results in an unstable or a qualitatively wrong solution. 

It is demonstrated in \cite{Maboudi:2016,Lall:2003iy,Carlberg:2014ky,Peng:2014di}, that if the basis for the reduced space is not chosen carefully, the symplectic symmetry of Hamiltonian and Lagrangian systems will be destroyed by model reduction. To resolve this issue, in \cite{Maboudi:2016,Peng:2014di,Lall:2003iy} a reduced order configuration space is constructed that inherits symmetries of the full configuration space. By using a proper time integrator scheme, the symmetries are preserved in the reduced system. A greedy-type algorithem is developed in \cite{Maboudi:2016} for construction of a basis for such a reduced configuration space.

Most models in engineering appear as a dissipative perturbation of a Hamiltonian system. In these systems, conservation of energy is taken as a fundamental principle of the system dynamics, while dissipative forces, e.g. friction, can change the energy of the system \cite{vanderSchaft:2014:PST:2693645.2693646}. As the energy is no longer preserved for such systems, existing methods can no longer be applied directly \cite{peng2016geometric}.

For dissipative and forced Hamiltonian system, Peng et al. \cite{peng2016geometric} suggest a symplectic model reduction method that preserves the Hamiltonian and the dissipative structure of the original system. However, since this method uses a symplectic integrator for a non-conservative system, there is no guarantee that the evolution of the energy is translated correctly to the reduced system.

In the context of network modeling and circuit simulation, considerable work has been done in the development of structure preserving, and in particular energy preserving, model reduction techniques. Model reduction for port-Hamiltonian systems are given in \cite{Polyuga:2010gj,beattie2011structure,chaturantabut2016structure} and the references therein. These methods use a Krylov or a Proper Orthogonal Decomposition (POD) approach to construct a reduced port-Hamiltonian system that preserves the passivity, and, thus, the stability of the original system. However, these methods do not generally guarantee the correct distribution of the energy among the energy consuming and energy storing units. Furthermore, over long time integration, accumulation of local errors might produce an erroneous solution.

In this paper, we present the Reduced Dissipative Hamiltonian (RDH) method as a structure preserving model reduction approach for dissipative Hamiltonian systems. A key difference between this method and the other existing methods is that the RDH enables the reduced system to be integrated using a symplectic integrator. By considering a canonical heat bath, also known as hidden strings \cite{Figotin:2006jy,Figotin:2005}, the reduced system is extended to a closed and conservative system. Therefore, a symplectic time integrator can be used to guarantee conservation of the system energy and the correct dissipation of energy. Furthermore, the hidden strings assure that the local errors in the dissipation of energy do not accumulate, resulting in a correct evolution of the system energy.  

What remains of this paper is organized as follow. Section \ref{sec:2} covers the required background on Hamiltonian systems, dissipative Hamiltonian systems, and the Hamiltonian extension. In Section \ref{sec:3} we discussion the symplectic model reduction for Hamiltonian systems and introduce the Reduced Dissipative Hamiltonian method. Accuracy, stability and efficiency of the RDH method is discussed in Section \ref{sec:4}, and illustrated through simulation of the dissipative wave equation and a linear port-Hamiltonian system of an electrical circuit. We offer conclusive remarks in Section \ref{sec:5}.

\section{Dissipative Hamiltonian Systems} \label{sec:2}
We first summarize needed concepts around the geometry of a symplectic linear vector space and then discuss Hamiltonian systems subject to dissipative forces. Finally, we introduce the Hamiltonian extension of dissipative Hamiltonian systems.

\subsection{Hamiltonian Systems} \label{sec:2.1} Suppose that $(\mathbb{R}^{2n},\Omega)$ is a symplectic linear vector space, where $\mathbb{R}^{2n}$ is a configuration space and $\Omega:\mathbb{R}^{2n}\times \mathbb{R}^{2n} \to \mathbb R$ is a closed, skew-symmetric and non-degenerate 2-form on $\mathbb{R}^{2n}$. Given a smooth Hamiltonian function $H:\mathbb{R}^{2n}\to \mathbb R$, \emph{Hamilton}'s equations of evolution are given as
\begin{equation} \label{eq:2.1}
	\begin{aligned}
	\dot {z}(t) &= \mathbb J_{2n} \nabla_{z} H, \\
	z(0) &= z_0,
	\end{aligned}
\end{equation}
where $z \in\mathbb R^{2n}$ is the state coordinates and $\mathbb J_{2n}$ is a $2n\times 2n$ matrix such that $\Omega(x,y) = x^T \mathbb J_{2n} y$, for all $x,y\in \mathbb R^{2n}$ \cite{Marsden:2010:IMS:1965128}. By using the \emph{Symplectic Gram-Schmidt} \cite{de2006symplectic} method, one can construct a coordinate system in which $\mathbb J_{2n}$ takes the form
\begin{equation} \label{eq:2.2}
	\mathbb{J}_{2n} = 
	\begin{pmatrix}
		0_n & I_n \\
		-I_n & 0_n
	\end{pmatrix},
\end{equation}
where $I_n$ and $0_n$ are the identity matrix and the zero matrix of size $n$, respectively. A main feature of Hamiltonian systems is the conservation of the Hamiltonian along the integral curves.
\begin{theorem} \label{theorem:2.1}
\cite{Marsden:2010:IMS:1965128} Consider the flow $\phi_t:\mathbb R\times \mathbb R^{2n} \to \mathbb R^{2n}$ of the Hamiltonian system (\ref{eq:2.1}). Then $H\circ \phi_t = H$.
\end{theorem}

In many physical problems, $H$ represents the system energy and is bounded from below. Here, to avoid difficulties with well-posedness of Hamiltonian systems,  we often assume that $H$ is a quadratic Hamiltonian, i.e., it takes the form $H(z) = \frac 1 2 z^T K^T K z$, where $K$ is a full rank $2n\times 2n$ matrix. This assumption leads to a linear system of evolution (\ref{eq:2.1}). We emphasise that the Hamiltonian extension in Section \ref{sec:2.2} and the Reduced Dissipative Hamiltonian method in Section \ref{sec:3.3} can be naturally extended to Hamiltonians of the form $H(z) = \frac 1 2 z^T K^T K z + g(z)$, where $g:\mathbb R^{2n} \to \mathbb R$ is an arbitrary function of $z$. 

Under general coordinate transformations, the equations of evolution may not take the form in (\ref{eq:2.1}). Let $(\mathbb R^{2n},\Omega)$ and $(\mathbb R^{2k},\Lambda)$ be two symplectic linear vector spaces. A linear transformation $\alpha :\mathbb R^{2n} \to \mathbb R^{2k}$ is called a \emph{symplectic transformation} \cite{Marsden:2010:IMS:1965128} if
\begin{equation}
	\Omega(x,y) = \Lambda(\alpha(x),\alpha(y)), \quad \text{for all } x,y\in \mathbb R^{2n}.
\end{equation}
In matrix notation, $A\in \mathbb R^{2n\times 2k}$ is called a \emph{symplectic matrix} if
\begin{equation}
	A^T \mathbb{J}_{2n} A = \mathbb{J}_{2k},
\end{equation}
where the superscript $T$ represents the transpose operator. As the symplectic 2-form is preserved under a symplectic transformation, the form of the equations of evolution remains invariant through a symplectic coordinate transformation \cite{Marsden:2010:IMS:1965128}. The \emph{symplectic inverse} of $A$, is a pseudo-inverse given by
\begin{equation}
	A^+ = \mathbb{\mathbb J}_{2k}^T A^T \mathbb J_{2n}.
\end{equation}
It is shown in \cite{Peng:2014di} that $A^+A = I_{2k}$ and that $(A^+)^T$ is a symplectic matrix. Furthermore, one easily checks that $AA^+$ is idempotent, i.e. it is a projection operator onto the column span of $A$. 

It is known that symplectic matrices are usually ill-conditioned \cite{Karow:2006cf}. Under some conditions on a symplectic space \cite{da2003introduction}, one can construct a symplectic basis which is also ortho-normal, and thus norm bounded. A basis which is both symplectic and orthogonal is called a \emph{ortho-symplectic basis}. We refer the reader to \cite{da2003introduction} for conditions on existence and construction of an ortho-symplectic basis.

It is natural to expect the numerical integrator that solves (\ref{eq:2.1}) also satisfy the conservation law in Theorem \ref{theorem:2.1}. However, common numerical integrators, e.g. the Runge-Kutta method, do not generally preserve the Hamiltonian. Symplectic numerical integrators are a class of numerical integrators for Hamiltonian systems that preserves the symplectic structure and ensure stability during long-time integration. The Str\"omer-Verlet time-stepping scheme
\begin{equation} \label{eq:2.3}
\begin{aligned}
	p_{n+1/2} &= p_n - \frac{\Delta t}{2} \nabla_qH(q_{n},p_{n+1/2}), \\
	q_{n+1} &= q_n + \frac{\Delta t}{2} \left( \nabla_pH(q_{n},p_{n+1/2}) + \nabla_pH(q_{n+1},p_{n+1/2}) \right),\\
	p_{n+1} &= p_{n+1/2} - \frac{\Delta t}{2} \nabla_qH(q_{n+1},p_{n+1/2}),
\end{aligned}
\end{equation}
where $p$ and $q$ are the canonical coordinates $z = (q^T,p^T)^T$, is an example of such numerical integrators. More information on the construction and application of symplectic integrators can be found in \cite{Hairer:1250576}.

\subsection{Dissipative Hamiltonian Systems and Hamiltonian Extensions} \label{sec:2.2}

Many systems in engineering and science appear as a perturbation of a Hamiltonian system, where the perturbation can be regarded as dissipation. In these systems, the energy tends to decrease over time, and thus, the conservation law in Theorem \ref{theorem:2.1} does not hold. Therefore, it is common to take the conservation of energy as a fundamental principle and consider the dissipative system coupled with a heat bath that absorbs the dissipated energy of the original system. 

To account for dissipation in a quadratic Hamiltonian $H(z) = \frac 1 2 z^T K^T K  z$, we rewrite (\ref{eq:2.1}) as a time dispersive and dissipative (TDD) \cite{Figotin:2006jy} system 
\begin{equation} \label{eq:2.4}
	\begin{aligned}
		& \dot {z} = \mathbb J_{2n} K^T f(t), \\
		& z(0) = z_0,
	\end{aligned}
\end{equation}
where $f$ is the solution to the Volterra integral equation \cite{corduneanu1991integral}
\begin{equation} \label{eq:2.4.1}
	f(t) + \int_0^t \chi(t-s) \cdot f(s)\ ds = K z.
\end{equation}
Here $\chi:\mathbb R^+\to \mathbb R^{2n\times 2n}$ is a bounded matrix valued function with respect to the Frobenius norm and is called the \emph{general susceptibility}. Note that the integral term in (\ref{eq:2.4.1}) accounts to the accumulation of the dissipation, whereas $\chi(s) = 0$ implies (\ref{eq:2.4}) is equivalent to (\ref{eq:2.1}). Furthermore, under suitable assumptions on $K$, both (\ref{eq:2.1}) and (\ref{eq:2.4}) are well-posed \cite{Figotin:2006jy}.

\begin{example} \label{example:2.1}
Consider the dynamics of the damped harmonic oscillator
\begin{equation} \label{eq:2.5}
	\ddot q + r \dot q + k q = 0
\end{equation}
where $k$ is the Hooke's constant and $r$ is the spring's damping factor. Note that without a damping term, (\ref{eq:2.5}) is a Hamiltonian system. The TDD formulation for the damped harmonic oscillator takes the form
\begin{equation} \label{eq:2.6}
	\dot q(t) = f(t), \quad \dot p(t) = - k q(t), \quad f(t) + \int_0^t r f(s) \ ds = p(t).
\end{equation}
Here $(q,p)$ are the canonical coordinates and the susceptibility is the constant function $r$.
\end{example}

It is shown in \cite{Figotin:2006jy,Figotin:2005} that under natural assumptions on the linear susceptibility $\chi(t)$ (see below), one can couple a TDD system of the form (\ref{eq:2.4}) with a canonical heat bath where the dissipated energy is captured in the heat bath in a canonical sense. In other words, one can construct a Hilbert space $\mathcal H$ and an isometric injection $I:\mathbb R^{2n} \to \mathbb R^{2n}\times \mathcal H^{2n}$ where the solution $z$ to (\ref{eq:2.4}) is the projection of $x$ onto $\mathbb R^{2n}$, and $x$ is the solution to
\begin{equation} \label{eq:2.7}
	\dot x = \mathcal J_{2n} \frac{\delta H_{\text{ex}}}{\delta x}.
\end{equation}
Here $H_{\text{ex}}:\mathbb R^{2n}\times \mathcal H^{2n} \to \mathbb R$ is an extended quadratic Hamiltonian function and $\mathcal J_{2n}$ is the symplectic operator defined on $\mathbb R^{2n}\times \mathcal H^{2n}$ respectively.

\begin{theorem}
Suppose that $K$ is full rank and $\chi(t)$ is symmetric. Then there is a quadratic extension to (\ref{eq:2.4}) of the form (\ref{eq:2.7}), if
\begin{equation} \label{eq:2.8}
	\text{Im}(\xi\hat{\chi}(\xi)) \geq 0, \quad \forall \xi = \omega + i\eta, \ \eta \geq 0,
\end{equation}
where $\hat{\chi}$ is the Fourier-Laplace transform of $\chi$
\begin{equation} \label{eq:2.9}
	\hat{\chi}(\xi) = \int_0^\infty e^{i\xi t} \chi(t)\ dt.
\end{equation}
\end{theorem}
\emph{Proof.} Here we prove the theorem for the case where $\chi$ is a constant symmetric matrix, where condition (\ref{eq:2.8}) corresponds to $\chi$ being positive semi-definite. We refer the reader to \cite{Figotin:2006jy} for the proof of the general case. Consider the Hamiltonian system
\begin{subequations}
\begin{align}
		\label{eq:2.10.a} & \dot{z}(t) = \mathbb J_{2n} K^T f(t), \\
		\label{eq:2.10.c} & \partial_t \phi(t,x) = \theta(t,x), \\
		\label{eq:2.10.b} & \partial_t \theta(t,x) = \partial_x^2 \phi(t,x) + \sqrt 2 \delta_0(x) \cdot \sqrt{\chi}  f(t), 
\end{align}
\end{subequations}
together with the initial condition
\begin{equation} \label{eq:2.10.1}
	z(0) = z_0,\quad \phi(0,\cdot) = 0, \quad \theta(0,\cdot) = 0.
\end{equation}
Here $\theta$ and $\phi$ are vector valued functions in $\mathcal H^{2n}$, $\delta_0(s)$ is the Dirac's delta function, $\sqrt{ \chi}$ is the matrix square root of $\chi$ and $f$ is the solution to the equation
\begin{equation} \label{eq:2.11}
	f(t) + \sqrt{2} \cdot \sqrt{ \chi } \phi(t,0) = Kz(t).
\end{equation}
To show that the Hamiltonian system (\ref{eq:2.10.a})-(\ref{eq:2.10.c}) is an extension to (\ref{eq:2.4}) in the sense discussed above, it is enough to show that the solution $f$ to equation (\ref{eq:2.11}) also satisfies (\ref{eq:2.4.1}). Equations (\ref{eq:2.10.c}) and (\ref{eq:2.10.b}) are equations for a vibrating string, and can be solved analytically
\begin{equation} \label{eq:2.12}
	\phi(t,x) = \frac {\sqrt 2} 2 \int_0^{t-|x|} \sqrt{\chi} f(s)\ ds,\quad \theta(t,x) = \frac{\sqrt 2}{2} \cdot \sqrt{\chi} f(t - |x|).
\end{equation}
Then, we recover (\ref{eq:2.4.1}) by substituting (\ref{eq:2.12}) into (\ref{eq:2.11}). The extended Hamiltonian $H_\text{ex}$ for the system (\ref{eq:2.10.a})-(\ref{eq:2.10.c}) takes the quadratic from
\begin{equation} \label{eq:2.13}
	H_\text{ex}(z,\phi,\theta) = \frac 1 2 \left( \| Kz - \phi(t,0) \|_2^2 + \| \theta(t) \|^2_{\mathcal H^{2n} } + \| \partial_x\phi(t)\|^2_{\mathcal H^{2n} }\right)
\end{equation}
where $\| \cdot \|_2$ is the Euclidean norm on $\mathbb R^{2n}$ and $\| \cdot \|_{\mathcal H^{2n}}$ is the induced norm from the inner product on $\mathcal H^{2n}$. 

Equations (\ref{eq:2.10.c}) and (\ref{eq:2.10.b}) are called the \emph{hidden strings}. The dissipation of energy in the original system (\ref{eq:2.4}) is carried away, as vibrations, along the added strings making the extended system conservative. The Hamiltonian extension of the damped harmonic oscillator in Example \ref{example:2.1} is exactly the Lamb model \cite{lamb:1900} which is a harmonic oscillator coupled with a vibrating string, and the tension in the string causes linear dissipation in the dynamics of the harmonic oscillator.

Note that the time integration of (\ref{eq:2.10.a})-(\ref{eq:2.10.c}) involves the integration of $f$ in (\ref{eq:2.12}). In general, the history of $f(t)$ must be stored and may cause storage limitation in long-time integration. However, we are interested solely in finding $z(t)$ which depends on $f$ at time $t$, and $\phi(t,0)$, i.e. the integral of the history of $f$. So by carefully choosing a quadrature rule that uses the same quadrature nodes as the time integrator we can avoid storing the history of $f$. For example for the trapezoidal rule, we recover the recursive relation
\begin{equation} \label{eq:2.14}
	\int_{0}^{t_n} f(s) \ ds \approx \frac{\Delta t}{2} f(t_n) + \frac{\Delta t}{2} f(t_{n-1}) + \int_{0}^{t_{n-1}} f(s) \ ds,
\end{equation}
where $\Delta t$ is the time step. The recursive relation in (\ref{eq:2.14}) suggests that storing the value of the integral term together with the state of $f$ in the previous time step suffices to evaluate the integral for the new time step. For other interpolation based quadrature rules, we can construct similar recursive rules of the form
\begin{equation}
	\int_{0}^{t_n} f(s) \ ds \approx \sum_{i=0}^{k} \omega_i f(t_{n-i})  + \int_{0}^{t_{n-k}} f(s) \ ds
\end{equation}
for some quadrature weights $\omega_i$, $i=1,\dots,k$ with $k\ll n$. Thus, time integration of (\ref{eq:2.10.a})-(\ref{eq:2.10.c}), only requires storage of $k$ evaluations of $f$.

\section{Model Order Reduction} \label{sec:3}
In this section we first explain the main results of \cite{Maboudi:2016,Peng:2014di} regarding model reduction of Hamiltonian systems and, subsequently we introduce the Reduced Dissipative Hamiltonian method.

\subsection{Symplectic Model Order Reduction} \label{sec:3.1}
Consider a Hamiltonian system of the form (\ref{eq:2.1}) together with a quadratic Hamiltonian of the form $H(z) = \frac 1 2 z^T K^T K z$. In this paper we focus on reducing the complexity of the numerical evaluation of (\ref{eq:2.1}) with respect to time $t$. Nevertheless, one can extend the results of this paper to a Hamiltonian system that depends on a set of physical or geometrical parameters belonging to a compact subset of a Euclidean space.

The main idea behind model order reduction is that the solution manifold, $\mathcal M = \{ z(t) : t \in [0,T] \}$ can be approximated by a low dimensional linear subspace \cite{hesthaven2015certified,quarteroni2015reduced}. A basis for such a subspace is called a \emph{reduced basis}, and its span is referred to as the \emph{reduced space} \cite{hesthaven2015certified}.

Suppose that a reduced symplectic basis $A \in \mathbb R^{2n\times 2k}$ is provided with $k \ll n$. The approximate solution to (\ref{eq:2.1}) in this basis is expressed as
\begin{equation} \label{eq:3.1}
	z =Ay,
\end{equation}
where $y\in \mathbb R^{2k}$ are the coordinates of the approximation with respect to the basis $A$. Substituting (\ref{eq:3.1}) into (\ref{eq:2.1}) yields
\begin{equation} \label{eq:3.2}
	A \dot y = \mathbb{J}_{2n} \nabla_{z} H(A y).
\end{equation}
Multiplying both sides with the symplectic inverse of $A$ and using the chain rule we obtain
\begin{equation} \label{eq:3.3}
	\dot y = A^+ \mathbb{J}_{2n} (A^+)^T \nabla_{y} H(A y).
\end{equation}
As $(A^+)^T$ is a symplectic matrix it implies that $A^+ \mathbb{J}_{2n} (A^+)^T = \mathbb{J}_{2k}$. By defining the reduced Hamiltonian $\tilde H : \mathbb{R}^{2k} \to \mathbb R$, as $\tilde H (y) = H(Ay)$, we recover the reduced system
\begin{equation} \label{eq:3.3.1}
	\begin{aligned}
	\dot {y}(t) &= \mathbb J_{2k} \nabla_{y} \tilde H, \\
	y(0) &= A^+ z_0,
	\end{aligned}
\end{equation}
Equation (\ref{eq:3.3.1}) is called the \emph{symplectic Galerkin projection} \cite{Peng:2014di} of the Hamiltonian system (\ref{eq:2.1}). Conventional model reduction techniques, e.g. Galerkin or Petrov-Galerkin methods \cite{hesthaven2015certified,quarteroni2015reduced}, do not yield a Hamiltonian reduced system and the reduced system does not necessarily preserve the conservation law in Theorem \ref{theorem:2.1} which results in a qualitatively wrong and often unstable solution \cite{Peng:2014di}. On the other hand the reduced system, obtained by the symplectic Galerkin projection, is a Hamiltonian system, and the system energy is therefore preserved over time \cite{Peng:2014di}. The following theorem guarantees the boundedness of the solution of the reduced system obtained by the symplectic Galerkin projection for quadratic Hamiltonians.

\begin{theorem}
\cite{Peng:2014di} Let $S$ be a bounded open subset of $\mathbb R^{2n}$ that contain $z_0$. Furthermore, assume that $H(z_0)<H(z)$ or $H(z_0)>H(z)$ for all $z\in \partial S$, the boundary of $S$. For a given symplectic basis $A$, provided $z_0$ is in the range of $A$, then both the original system and the reduced system obtained by the symplectic Galerkin projection are bounded.
\end{theorem}

In the next section we introduced the greedy method to construct a symplectic reduced basis \cite{Maboudi:2016}.

\subsection{The Greedy Approach to Symplectic Basis Generation} \label{sec:3.2}
Greedy generation of the reduced basis is an iterative procedure which, in each iteration, adds the two best possible basis vectors to the symplectic basis to enhance overall accuracy. Suppose that at the $k$-th generic step of the greedy basis selection, an ortho-symplectic basis $A_{2k} = \{ e_1,\dots,e_k,f_1,\dots,f_k\}$ of size $2k$ is provided. The first step in an iteration of the greedy basis selection comprises finding $t^{k+1}\in[0,T]$ such $z(t^{k+1})$ is worst approximated by the current reduced space. In other words
\begin{equation} \label{eq:3.4}
	t^{k+1} := \underset{t\in [0,T]}{\text{argmax }} \| z(t) - A_{2k}{A_{2k}}^+z(t) \|.
\end{equation}
Suppose that $e_{k+1}$ is the vector obtained by $\mathbb J_{2n}$-orthogonalization (the symplectic Gram-Schmidt process \cite{Salam2014}) of $z(t^{k+1})$ with respect to $A_{2k}$. Then the enriched basis $A_{2k+2}$ takes the form
\begin{equation}
	A_{2k+2} = \{ e_1,\dots,e_k, e_{k+1},f_1,\dots,f_k,\mathbb J_{2n}^T e_{k+1}\}.
\end{equation}
It is easily checked that $A_{2k+2}$ is ortho-symplectic. Furthermore, it is shown in \cite{Maboudi:2016} that under natural assumptions on the solution manifold $\mathcal M$, the greedy method converges exponentially fast.

For parametric problems, evaluation of the projection error in (\ref{eq:3.4}) for the entire parameter space is computationally demanding. One can use the loss in the Hamiltonian function as a cheap surrogate to the projection error, i.e.
\begin{equation}
	\omega^{k+1} := \underset{\omega\in \Omega}{\text{argmax }} | H(z(\omega)) - H(AA^+z(\omega)) |.
\end{equation}
Here $\Omega$ is a closed and bounded set of parameters for the original Hamiltonian system. Note that $\omega^{k+1}$ can be identified prior to time integration since the Hamiltonian function is constant in time, i.e. $H(z) = H(z_0)$ and $H(AA^+z) = H(AA^+z_0)$. We summarize the greedy method for symplectic basis generation in Algorithm \ref{alg:added1} and refer the reader to \cite{Maboudi:2016} for further details.

\begin{algorithm} 
\caption{The greedy algorithm for generation of a symplectic basis} \label{alg:added1}
{\bf Input:} Tolerated projection error $\delta$, initial condition $ z_0$
\begin{enumerate}
\item $t^1 \leftarrow t=0$
\item $e_1 \leftarrow z_0$
\item $A \leftarrow [e_1,\mathbb J^T_{2n}e_1]$
\item $k \leftarrow 1$
\item \textbf{while} $\| z(t) - A_{2k}{A_{2k}}^+z(t) \| > \delta$ for all $t \in [0,T]$
\item \hspace{0.5cm} $t^{k+1} := \underset{t\in [0,T]}{\text{argmax }} \| z(t) - A_{2k}{A_{2k}}^+z(t) \|$
\item \hspace{0.5cm} $\mathbb J_{2n}$-orthogonalize $ z(t^{k+1})$ to obtain $e_{k+1}$
\item \hspace{0.5cm} $A \leftarrow [e_1,\dots ,e_{k+1} , \mathbb J^T_{2n}e_1,\dots,\mathbb J^T_{2n}e_{k+1}]$
\item \hspace{0.5cm} $k \leftarrow k+1$
\item \textbf{end while}
\end{enumerate}
\vspace{0.5cm}
{\bf Output:} Symplectic basis $A$.
\end{algorithm}

\subsection{The Reduced Dissipative Hamiltonian Method} \label{sec:3.3}

Since the symplectic model reduction in Section \ref{sec:3.2} is based on the conservation law in Theorem \ref{theorem:2.1}, it can no longer be applied to dissipative Hamiltonian systems. Instead in the Reduced Dissipative Hamiltonian method, we considers a Hamiltonian extension to a dissipative Hamiltonian system to construct a closed system. A symplectic model reduction can then be naturally applied to conserve the total energy.

Consider a dissipative Hamiltonian system of the form (\ref{eq:2.4}) with a quadratic Hamiltonian, $H(z) = z^TK^TKz$. Since $K^TK$ is symmetric and positive definite, it has a unique Cholesky factorization $K^TK = L^T L$ where $L$ is upper triangular \cite{strang09}. So we can write 
\begin{equation} \label{eq:3.7}
H(z) = z^T L^T L z.
\end{equation}
Further, suppose that the solution $z(t)$ lies on a low-dimensional symplectic subspace such that $z = Ay$, where $A$ is an ortho-symplectic matrix of the size $2n\times 2k$ and $y$ is the expansion coefficients of $z$ in the basis of $A$. The Hamiltonian $H$ then takes the form
\begin{equation} \label{eq:3.8}
	H(z) = H(Ay) = y^T A^T L^T L A y.
\end{equation}
Since $A^T L^T L A$ is symmetric and positive definite, it too has a unique Cholesky factorization $\tilde L^T \tilde L$ where $\tilde L = A^T L A$ is an upper triangular matrix of size $2k \times 2k$. Writing (\ref{eq:2.4}) in terms of the reduced coordinates $y$ and the Hamiltonian (\ref{eq:3.7}) reads
\begin{equation} \label{eq:3.9}
		A\dot{y}(t) = \mathbb J_{2n} L^T f(t),
\end{equation}
together with the complementary equation
\begin{equation} \label{eq:3.10}
	f(t) + \sqrt 2 \cdot \sqrt{\chi} \phi(t,0) = LAy.
\end{equation}
Multiplying (\ref{eq:3.9}) with $A^+$ and (\ref{eq:3.10}) with $A^T$ yields
\begin{align} \label{eq:3.11}
	& \dot y(t) = \mathbb J_{2k} A^T L^T f(t), \\
	& A^T f(t) + \sqrt{2} A^T \sqrt{\chi} \phi(t,0) = A^T L A y,
\end{align}
where we use the fact that $A^+\mathbb J_{2n} = \mathbb{J}_{2k} A^T$. If we define $f = A \tilde f$, $\phi = A \tilde \phi$, $\theta = A\tilde \theta$ and the \emph{reduced susceptibility} as $\tilde \chi = A^T \chi A$ we recover the reduced Hamiltonian system
\begin{subequations}
\begin{align}
		\label{eq:3.12.a} & \dot{y}(t) = \mathbb J_{2k} {\tilde L}^T \tilde f(t), \\
		\label{eq:3.12.c} & \partial_t \tilde \phi(t,x) = \tilde \theta(t,x),\\
		\label{eq:3.12.b} & \partial_t \tilde \theta(t,x) = \partial_x^2 \tilde \phi(t,x) + \sqrt 2 \delta_0(x) \cdot \sqrt{\tilde \chi}  \tilde f(t),
\end{align}
\end{subequations}
together with the auxiliary equation
\begin{equation} \label{eq:3.13}
	\tilde f(t) + \sqrt{2} \sqrt{\tilde \chi} \tilde \phi(t,0) = \tilde L y.
\end{equation}
Equations (\ref{eq:3.12.a})-(\ref{eq:3.12.b}) is a Hamiltonian system on the symplectic linear vector space $\mathbb R^{2k} \times \mathcal H^{2k}$ and contributes to the \emph{reduced TDD system}
\begin{equation}
	\dot {y} = \mathbb J_{2k} \tilde L^T \tilde f(t), \quad \tilde f(t) + \int_0^t \tilde \chi\cdot \tilde f(s)\ ds = \tilde L y.
\end{equation}
Therefore, the system energy will be conserved along integral curves of (\ref{eq:3.12.a})-(\ref{eq:3.12.b}).

We point out that the transformation that connects (\ref{eq:2.10.a})-(\ref{eq:2.10.b}) to (\ref{eq:3.12.a})-(\ref{eq:3.12.b}) is given by
\begin{equation}
	\mathcal A = \begin{pmatrix}
		A& 0 \\
		0& A
	\end{pmatrix} : \mathbb R^{2n} \times \mathcal H^{2n} \to \mathbb R^{2k} \times \mathcal H^{2k}.
\end{equation}
This is a symplectic transformation, since $\mathcal A^T \mathcal J_{2n} \mathcal A = \mathcal J_{2k}$. Furthermore, the dissipation of energy in the reduced system only depends on the reduced susceptibility. Thus, the choice of $A$ should be independent of the hidden strings $(\phi, \theta)$. In other words, if the reduced space is chosen to be a symplectic subspace, then the actions of model reduction and Hamiltonian extension commute. We summarize the algorithm for model reduction of dissipative Hamiltonian systems in Algorithm \ref{alg:3.1}.

\begin{algorithm}
\caption{The Reduced Dissipative Hamiltonian Method (RDH)} \label{alg:3.1}
\begin{enumerate}
	\item Construct the Hamiltonian extension (\ref{eq:2.10.a})-(\ref{eq:2.10.b}) to the original TDD system (\ref{eq:2.4}).
	\item Collect the snapshots $z(t_i)$, $i=1,\dots,N$ through time integration of the extended Hamiltonian.
	\item Construct an ortho-symplectic basis $A$.
	\item Define $\tilde L = A^T L A$, $\tilde \chi = A^T \chi A$ and construct the reduced dissipative Hamiltonian system (\ref{eq:3.12.a})-(\ref{eq:3.12.b})
\end{enumerate}
\end{algorithm}

Note that Algorithm \ref{alg:3.1} does not depend on the choice of the method to construct an ortho-symplectic basis $A$. Thus, for basis generation, the greedy approach introduced in section \ref{sec:3.2} or an SVD-based method, e.g. the cotangent lift \cite{Peng:2014di}, can be applied.

The main advantage of the RDH method compared to the existing methods is that it enables the reduced system to be integrated using a symplectic integrator. The reduced system constructed using the RDH is a closed Hamiltonian system, therefore the conservation law in Theorem \ref{theorem:2.1} holds and a symplectic integrator guarantees that the total energy is preserved in the reduced system. Alternative methods, e.g. \cite{peng2016geometric,Polyuga:2010gj,beattie2011structure}, either integrate the reduced system with a non-symplectic integrator, or do not construct a closed reduced system which result in accumulation of local errors or unstable solution during long time integration, respectively \cite{Hairer:1250576}.

\section{Numerical Results} \label{sec:4}
In the following we illustrate the performance of the method through the reduced order model of the dissipative wave equation and a port-Hamiltonian model for a dissipative circuit.

\subsection{Dissipative wave equation} \label{sec:4.1}

Consider the dissipative linear wave equation
\begin{equation} \label{eq:4.1}
	\left\{
	\begin{aligned}
		q_{t}(t,x) &= p(t,x), \\
		p_{t}(t,x) &= c^2 q_{xx}(t,x) - r(x)  p(t,x) , \\
		q(0,x) &= q_0(x), \\
		p(0,x) &= 0.
	\end{aligned}
	\right.
\end{equation}
where $x$ belongs to a one-dimensional torus of length $L$ and $r:[0,1]\to[0,1]$ is a positive semi-definite real valued function. 

We discretize the torus into $N_{\Delta x}$ equidistant points and define $\Delta x = L/N_{\Delta x}$, $x_i = i\Delta x$, $q_i=q(t,x_i)$ and $p_i=p(t,x_i)$ for $i = 1, \dots, N_{\Delta x}$. The discretization of $r$ corresponds to a diagonal and semi-positive definite matrix $r_\Delta$. Furthermore, we discretize (\ref{eq:4.1}) using a standard central finite differences schemes to obtain
\begin{equation} \label{eq:4.2}
	\dot z = \mathbb J_{2n} K^T K z - R z,
\end{equation}
where $z = (q_1,\dots,q_{N_{\Delta x}},p_1,\dots,p_{N_{\Delta x}})$ and $K$ and $R$ are given as
\begin{equation} \label{eq:4.3}
	K^T K =
	\begin{pmatrix}
		I & 0 \\
		0 & c^2D_x^TD_x
	\end{pmatrix} , \quad
	R =
	\begin{pmatrix}
		0 & 0 \\
		0 & r_\Delta
	\end{pmatrix},
\end{equation}
with $D_x^TD_x = D_{xx}$ as the central finite differences matrix operator. Writing (\ref{eq:4.2}) in a TDD formulation yields
\begin{equation} \label{eq:4.4}
	\dot z = \mathbb J_{2n} K^T f(t), \quad f(t) + R \int_0^t f(s) \ ds = K z.
\end{equation}
Since $R$ is not time dependent, it commutes with the integration operator. The Hamiltonian extension of (\ref{eq:4.4}), then takes the form (\ref{eq:2.10.a})-(\ref{eq:2.10.b}).

The initial condition used is given by
\begin{equation} \label{eq:4.5}
	q_i(0) = h( 10\times|x_i - \frac{1}{2}| ), \quad p_i = 0, \quad i=1,\dots,N,
\end{equation}
where $h(s)$ is the cubic spline function
\begin{equation} \label{eq:4.6}
h(s) = 
\left\{
\begin{aligned}
& 1 - \frac{3}{2}s^2 + \frac{3}{4}s^3, \quad & 0\leq s \leq 1, \\
& \frac{1}{4}(2-s)^3, & 1< s \leq 2, \\
& 0, & s > 2.
\end{aligned}
\right.
\end{equation}
For the numerical time integration of the extended Hamiltonian system, the Str\"omer-Verlet time stepping scheme (\ref{eq:2.3}) is used. In each time step, the system of linear equations (\ref{eq:2.11}) is solved to recover $z$. System parameters are summarized below.
\vspace{0.5cm}
\begin{center}
\begin{tabular}{|l|l|}
\hline
Domain length & $L = 1$ \\
No. grid points & $N = 500$ \\
Space discretization size & $\Delta x = 0.002$ \\
Time discretization size & $\Delta t = 0.002$ \\
Wave speed & $c^2 = 0.1$ \\
\hline
\end{tabular}
\end{center}
\vspace{0.5cm}
The first numerical experiment corresponds to an inhomogeneous dissipative media. Here, $r_{\Delta}$ is a diagonal matrix with diagonal elements $r_i := 0.1 + 0.9(i/N_{\Delta x})$, for $i=1,\dots,N_{\Delta x}$.

\begin{figure}[t]
\begin{tabular}{cc}
\includegraphics[width=0.5\textwidth]{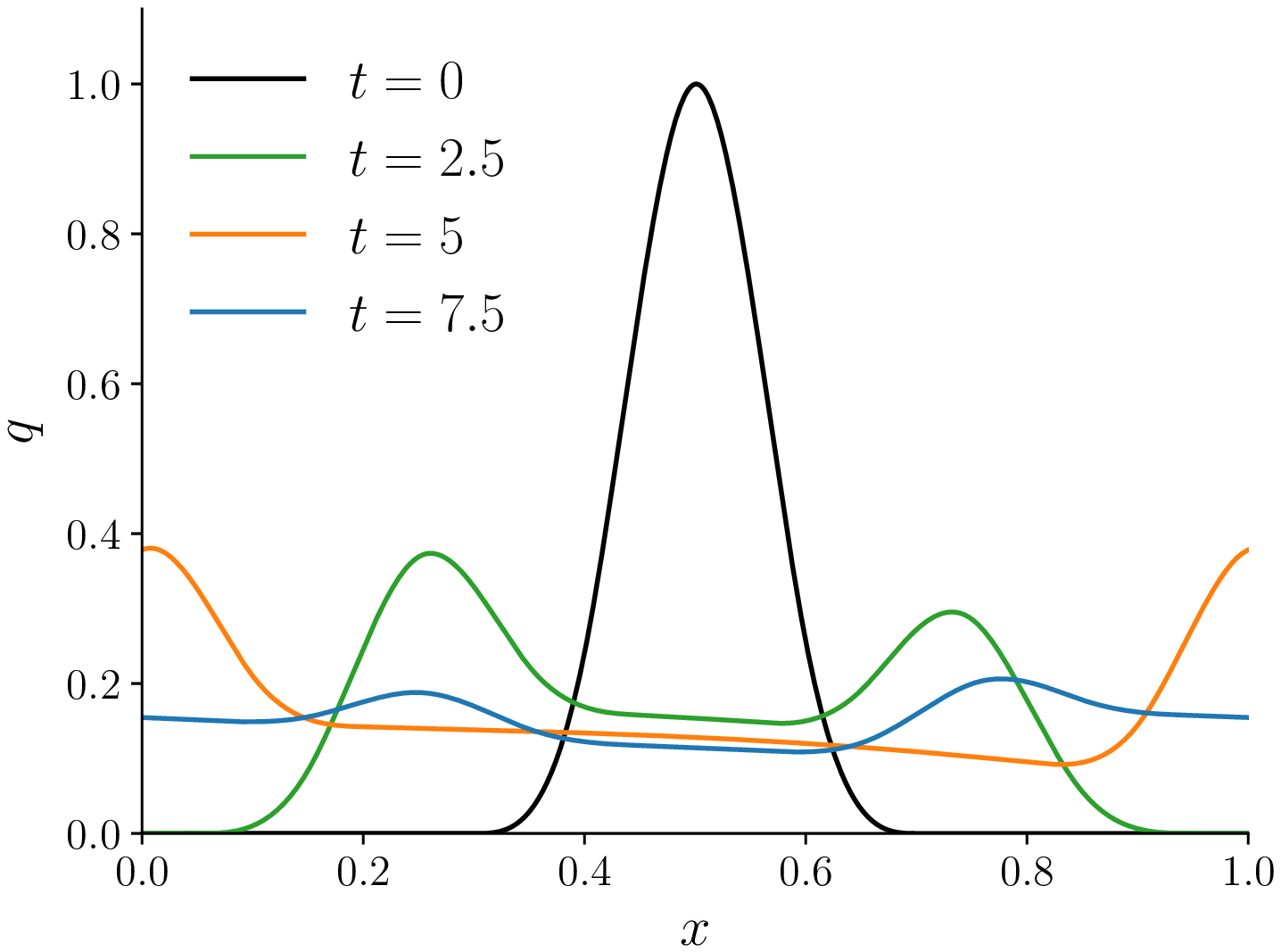} & 
\includegraphics[width=0.5\textwidth]{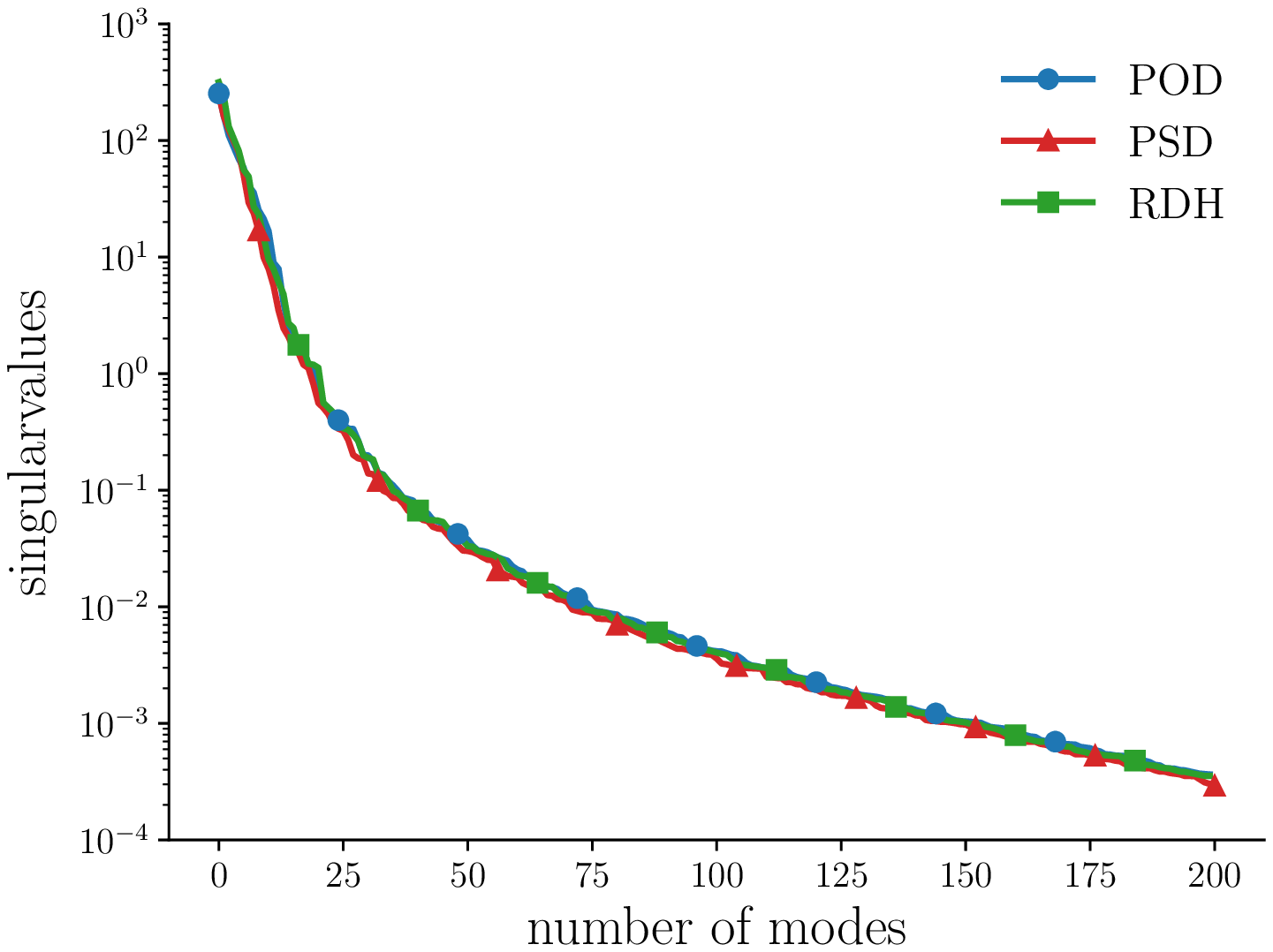} \\
(a) & (b) \\
\includegraphics[width=0.5\textwidth]{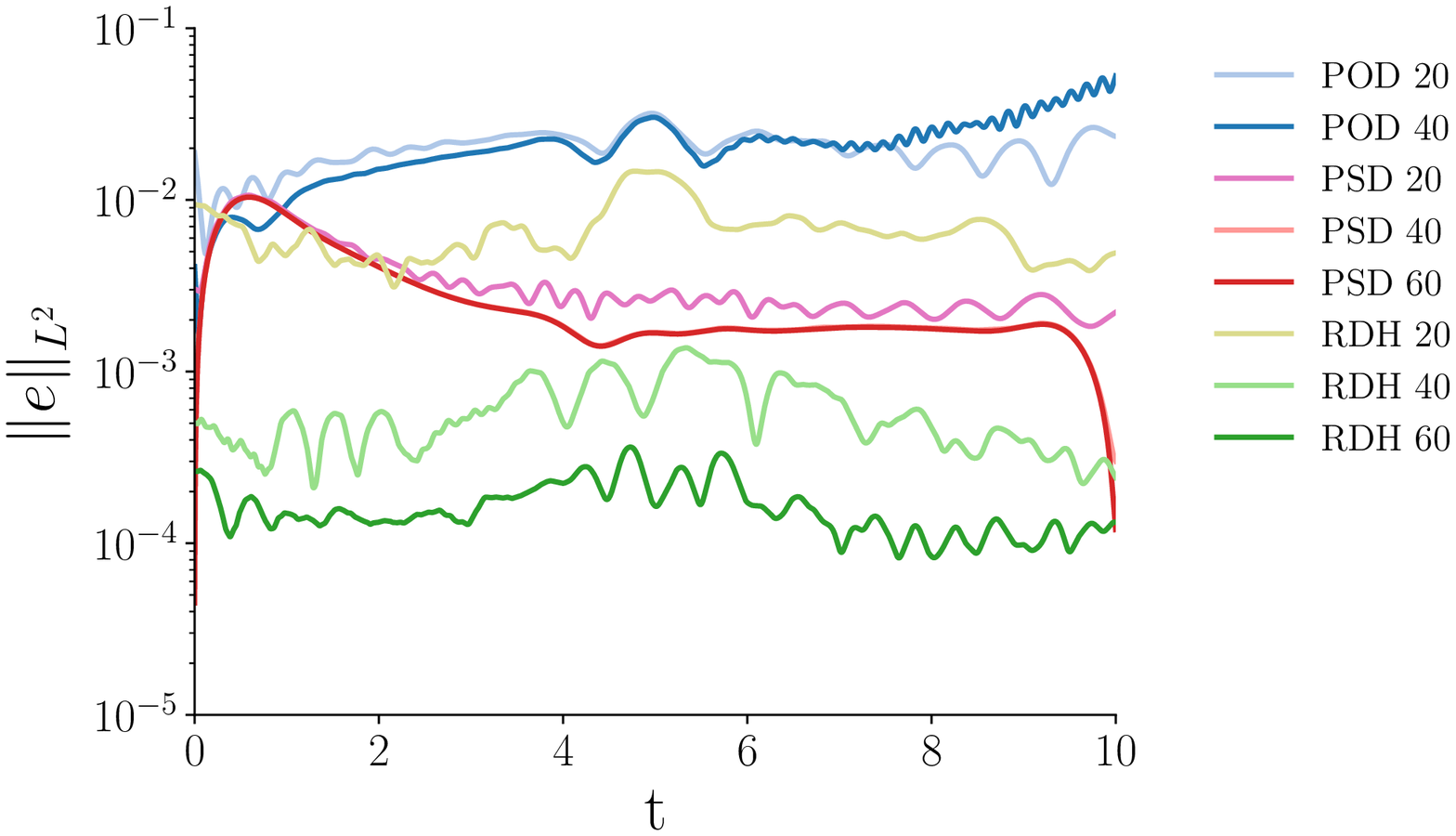} & 
\includegraphics[width=0.5\textwidth]{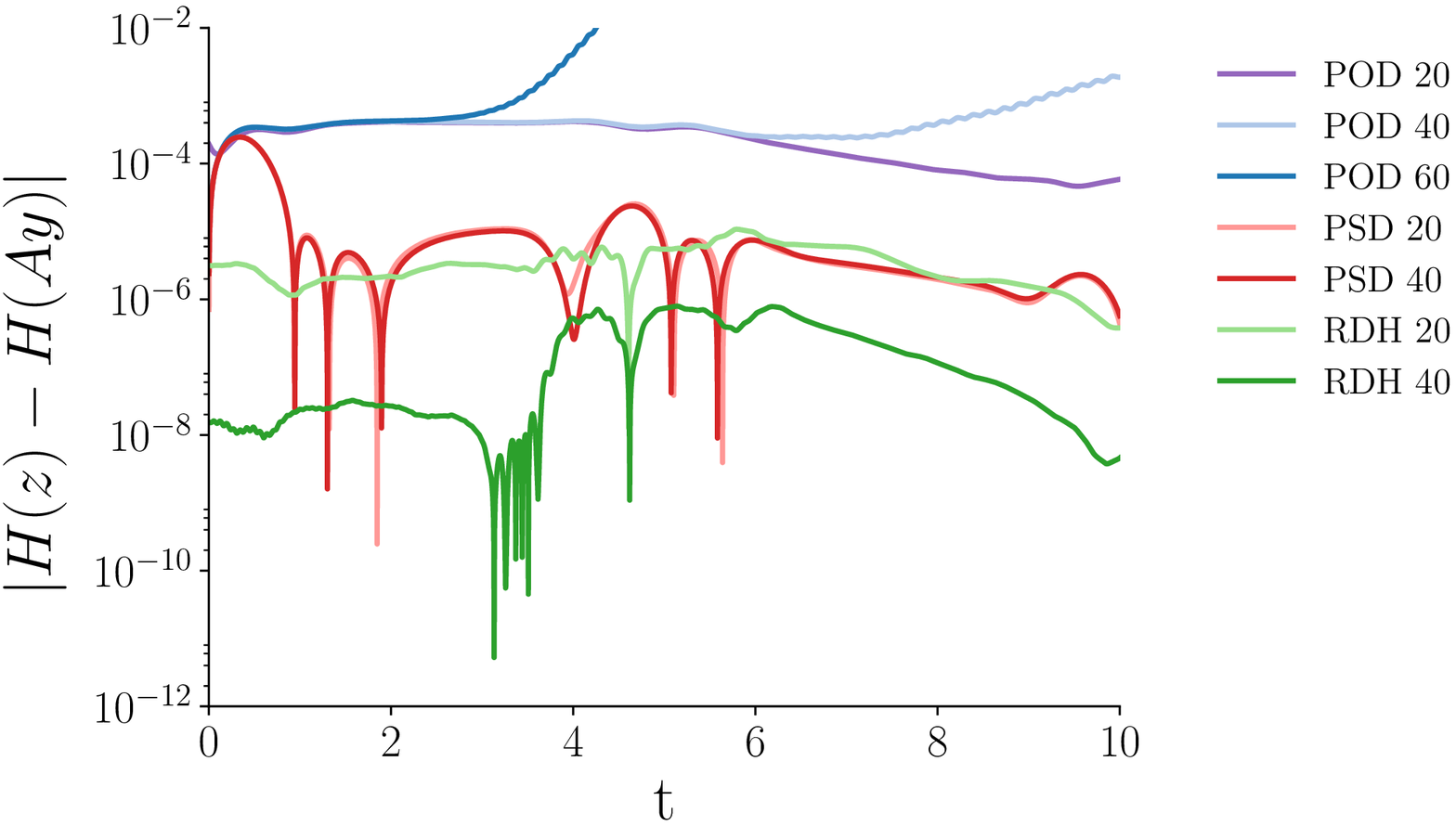} \\
(c) & (d) \\
\includegraphics[width=0.55\textwidth]{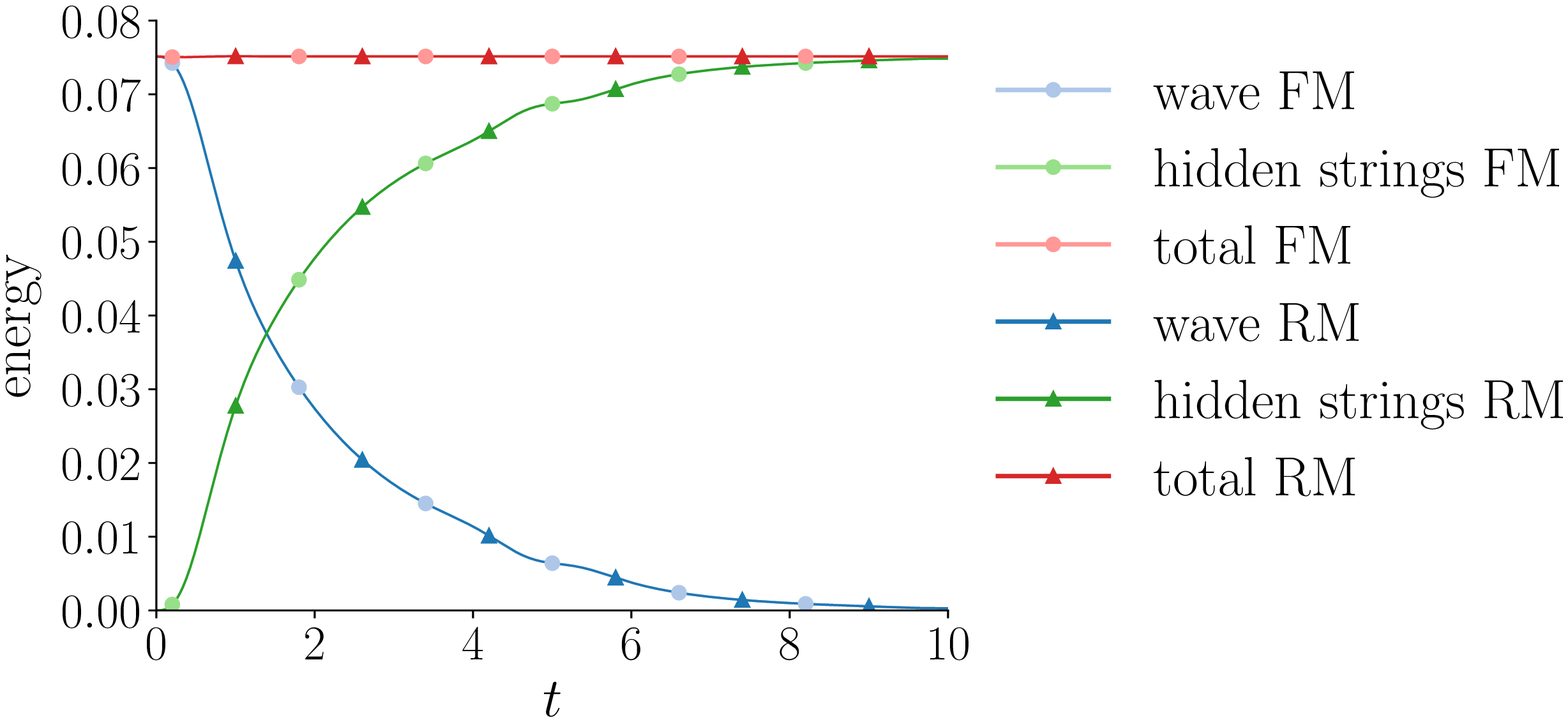} &
\includegraphics[width=0.45\textwidth]{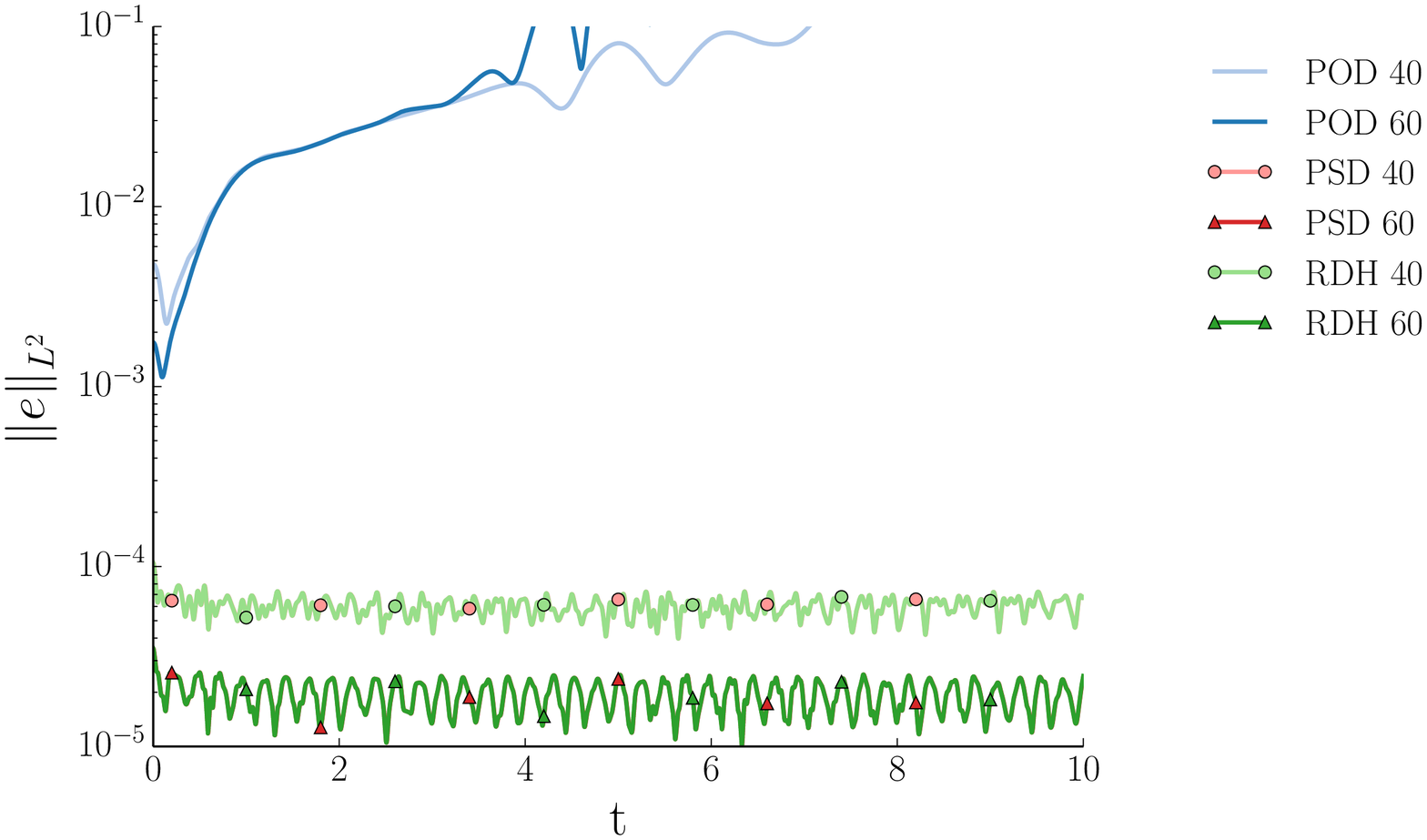} \\
(e) & (f) 
\end{tabular}
\caption{(a) The solution to the original dissipative wave equation (\ref{eq:4.1}), (b) The decay of the singular values for the POD, the PSD, and the RDH methods, (c) The $L^2$-error for the different methods, (d) Evolution of error in the Hamiltonian for different methods, (e) Energy preservation of the Hamiltonian extension for the original and the reduced system. ``FM'' and ``RM'' refer to the full model and the reduced model, respectively. (f) The $L^2$-error between the solution to the reduced system and the full system in a near-zero dissipation regime.} \label{fig:4.1}
\end{figure}

Figure \ref{fig:4.1}.(a) shows the solution of the original dissipative wave equation (\ref{eq:4.1}) at $t \in \{0,2.5,5,7.5\}$. For a nonzero $r_\Delta$ the solution will converge to $(q(t=\infty,x),p(t=\infty,x)) = (\rho,0)$ where $\rho$ is the center of mass of $q_0$. 

We construct the RDH reduced system according to the Algorithm \ref{alg:3.1}. The performance of the method is then compared to the POD and the method proposed in \cite{peng2016geometric}, referred to as the PSD.

Figure \ref{fig:4.1}.(b) illustrates the decay of the singular values of the snapshot matrix \cite{hesthaven2015certified}, for the POD, PSD, and the RDH methods. Note that the snapshots for the PSD and the RDH are different since they have different canonical representations. The fast decay of the eigenvalues in all methods is a strong indicator for the existence of a low dimensional reduced system. The reduced bases are then constructed using 20, 40 and 60 number of modes.

The $L^2$-error between the full system and the RDH, the PSD, and the POD methods are presented in Figure \ref{fig:4.1}.(c). We notice that the symplectic methods provide a more accurate solution when compared to the POD method. In fact, the POD method does not yield a stable reduced system.  Furthermore, it is seen that enriching the PSD reduced basis does not yield a significant improvement in the accuracy of the reduced system. This happens as the PSD method, numerically integrate a non-conservative system with a symplectic integrator. This results in an incorrect evolution of the energy and eventually, in a qualitatively wrong numerical solution.

On the other hand, we notice that the RDH method with 40 modes provides a significantly more accurate solution compared to the PSD method with 60 modes. The RDH method provides a conservative reduced system where the dissipated energy is absorbed by the hidden strings and the conservation of the energy is then guaranteed by using a symplectic integrator. Therefore, we observe remarkable increase in the accuracy by enriching the RDH reduced basis.

Figure \ref{fig:4.1}.(d) shows the conservation of the energy in the different methods. The conservation law expressed in Theorem \ref{theorem:2.1} is destroyed through the POD model reduction and as a consequence we observe blow-up of the system energy. The symplectic methods preserves the energy significantly better. As discussed above, enriching the PSD basis does not significantly improve the preservation of energy. On the contrary, the RDH provides a substantial improvement in the accuracy of the energy.

In Figure \ref{fig:4.1}.(e) we show the transfer of the energy from the TDD system to the hidden strings, for the full system and the RDH reduced system. We notice that the RDH method preserves the total energy of the extended Hamiltonian system. Furthermore, the transfer of energy to the hidden strings in the full model is correctly translated in the reduced system.

The second numerical experiment is the dissipative wave equation (\ref{eq:4.1}) in a near-zero dissipation regime. The numerical setting is taken to be identical to the previous numerical experiment, but with the difference that $r_i = 10^{-5}$, for $i=1,\dots,N_{\Delta x}$. Figure \ref{fig:4.1}.(f) shows the $L^2$-error between the solution to the reduced system and the full system, for the POD, the PSD, and the RDH methods. We notice that the POD does not yield a stable reduced system as the symplectic structure is lost via model reduction. Furthermore we notice that error for the PSD and the RDH coincide as the two methods become identical as $\| \chi \|_{\infty}\to 0$.

\subsection{The sine-Gordon equation}
Consider the one-dimensional dissipative nonlinear wave equation
\begin{equation} \label{eq:added.4.1}
	\left\{
	\begin{aligned}
		q_{t}(t,x) &= p(t,x), \\
		p_{t}(t,x) &= q_{xx}(t,x) - \sin(q) - r(x)  p(t,x), \\
		q(0,x) &= q_0(x), \quad q(t,0) = a, \quad q(t,L) = b,\\
		p(0,x) &= p_0(x).
	\end{aligned}
	\right.
\end{equation}
defined on a domain of length $L$, which is known as the sine-Gordon equation. In the absence of dissipation, $r(x) = 0$, the \emph{kink} solution to (\ref{eq:added.4.1}) is given as
\begin{equation} \label{eq:added.4.2}
	q(t,x) = 4\arctan\left( \exp \left( \frac{(x-x_0 - vt)}{\sqrt{1-v^2}} \right) \right),
\end{equation}
where $|v|< 1$ is the wave speed. In the presence of dissipation, where $r(x)\geq 0$, the traveling wave de-accelerates and stops. The TDD formulation for (\ref{eq:added.4.1}) takes the form
\begin{equation} \label{eq:added.4.3}
	\dot z = \mathbb J_{2n} K^T f(t) + \mathbb J_{2n} g(z) - \mathbb J_{2n} z_{\text{bd}}, \quad f(t) + R \int_0^t f(s) \ ds = K z.	
\end{equation}
where $z$, $K$, $R$ and $f$ are defined similar to (\ref{eq:4.4}), with $r_{\Delta} = rI_{n}$, and 
\begin{equation}
g(z) = (\sin(q_1),\dots,\sin(q_{N_{\Delta x}}),0,\dots,0)^T,
\end{equation}
and $z_{\text{bd}}$ is the term corresponding to the Dirichlet boundary condition. Note that the extended Hamiltonian $H_{ex}$ takes the form
\begin{equation}
	H_\text{ex}(z,\phi,\theta) = \frac 1 2 \left( \| Kz - \phi(t,0) \|_2^2 + \| G(z) \|_2^2 + \| \theta(t) \|^2_{\mathcal H^{2n} } + \| \partial_x\phi(t)\|^2_{\mathcal H^{2n} }\right),
\end{equation}
where $G(z)$ is a potential for $g(z)$ given as $G(q_i) = 1 - \sin(q_i)$, for $i=1,\dots,N_{\Delta_x}$. System parameters are summarized below

\vspace{0.5cm}
\begin{center}
\begin{tabular}{|l|l|}
\hline
Domain length & $L = 50$ \\
No. grid points & $N = 500$ \\
Space discretization size & $\Delta x = L/N$ \\
Time discretization size & $\Delta t = 0.02$ \\
Wave speed & $v = 0.5$ \\
Boundary conditions & $a = 0$, $b=1$\\
Dissipation coefficient & $r = 0.1$ \\
\hline
\end{tabular}
\end{center}
\vspace{0.5cm}

The RDH reduced system is constructed following Algorithm \ref{alg:3.1}. To reduce the complexity of the nonlinear term, we used the symplectic discrete empirical interpolation method (SDEIM) \cite{Maboudi:2016}. The performance of the method is then compared to the PSD and the POD where the SDEIM proposed in \cite{Peng:2014di} and the classical DEIM \cite{Chaturantabut:2010cz} is applied to reduce the complexity of the nonlinear term, respectively.

\begin{figure}[t]
\begin{tabular}{cc}
\includegraphics[width=0.5\textwidth]{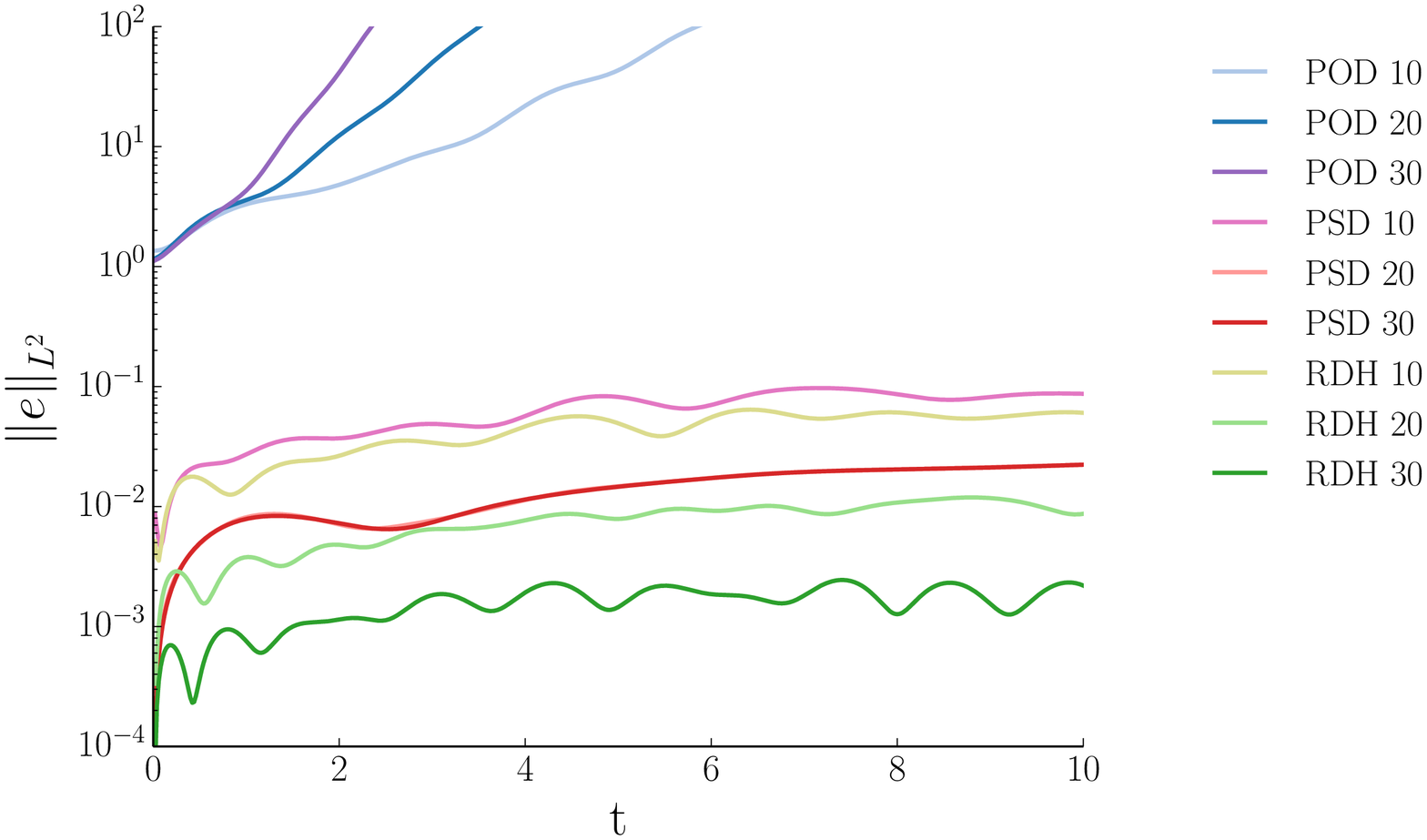} & 
\includegraphics[width=0.5\textwidth]{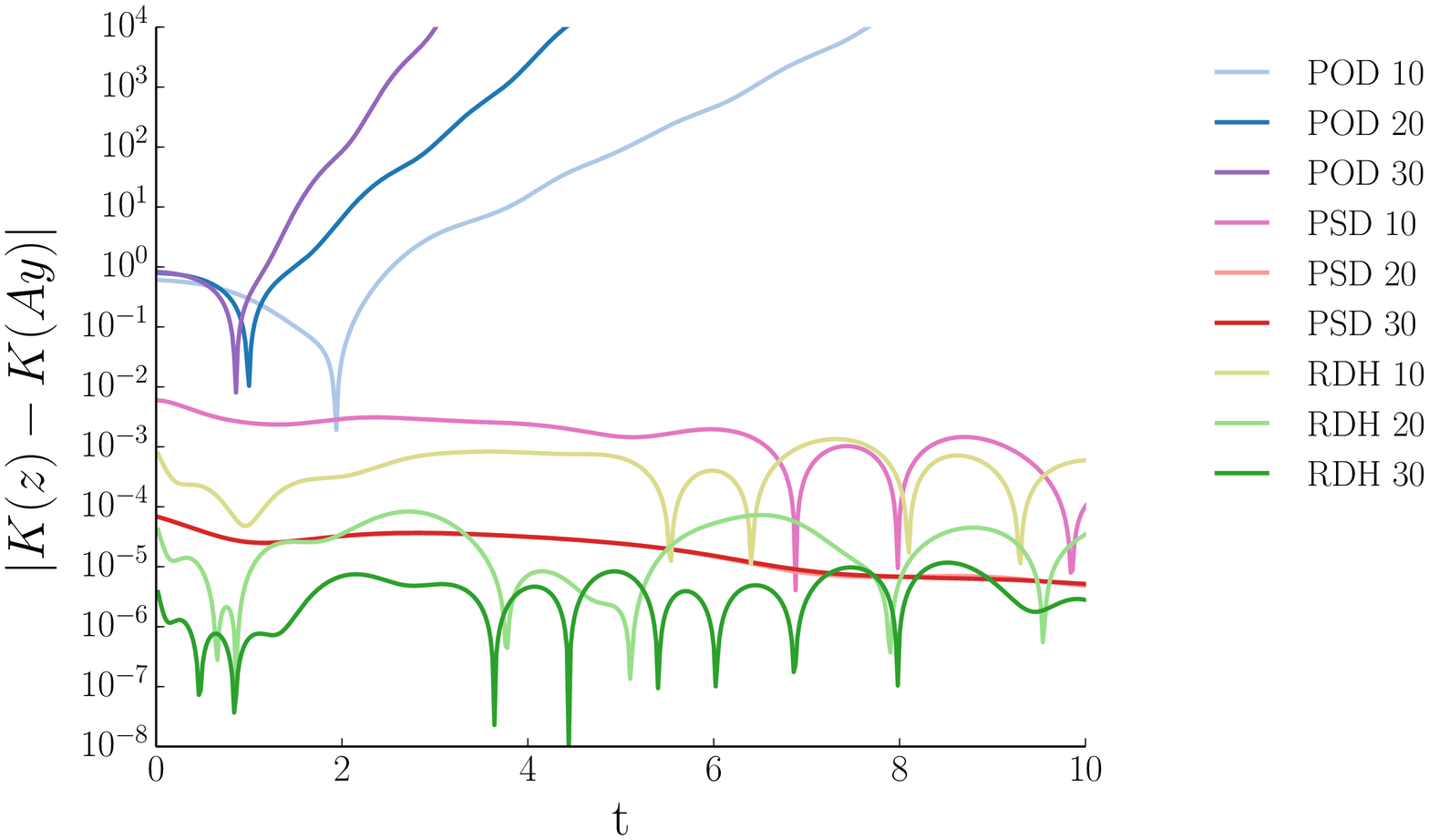} \\
(a) & (b)
\end{tabular}
\caption{(a) The $L^2$-error for the different methods, (a) Evolution of error in the kinetic energy for different methods.} \label{fig:added4.1}
\end{figure}

Figure \ref{fig:added4.1}.(a) shows the The $L^2$-error between the full system and the RDH, the PSD, and the POD methods. Although the Hamiltonian system of the sine-Gordon equation is nonlinear, the errors for the different methods show a similar behavior as those in Section (\ref{sec:4.1}). We observe that the POD does not yield a stable reduce system while the symplectic methods provide a high accuracy. Furthermore, we notice that enriching the PSD basis does not significantly enhance the accuracy of the method.

The evolution of error in the kinetic energy $K(p) = \|p\|_2^2/2$ is illustrated in Figure \ref{fig:added4.1}.(b). We see that the POD does not conserve the evolution of the kinetic energy. The RDH method conserves the kinetic energy of the system with a higher accuracy than the PSD method. Furthermore, the accuracy of the RDH method is better scaled under enrichment of the reduced basis, compared to the PSD method. 

It is observed in Figure 1 that the symplectic treatment of the nonlinear terms is essential in correct model reduction of Hamiltonian systems. In addition, the SDEIM can be combined with the RDH method to construct a reduced Hamiltonian system that can be integrated using a symplectic integrator. Thus, the combination preserves the system energy and the symplectic symmetry of Hamiltonian systems.

\subsection{Port-Hamiltonian Systems}
Port-Hamiltonian systems are popular in network modeling and electrical engineering. In the framework of port-Hamiltonian modelling, energy conservation and Hamiltonian structure is the fundamental principle of the dynamics of the system. Ports in the system network then allows the exchange of energy with the environment in the form of sources, capacitors, and dissipations \cite{vanderSchaft:2014:PST:2693645.2693646}. Port-Hamiltonian systems can be viewed as a forced and dissipative Hamiltonian system.
\begin{figure}[t]
\begin{center}
	\includegraphics[width=0.7\textwidth]{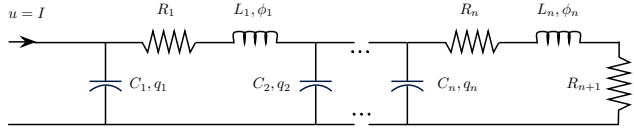}
\end{center}
\caption{$n$-dimensional ladder network} \label{fig:4.2}
\end{figure}

Consider the $n$-dimensional linear ladder network in Figure \ref{fig:4.2}. Here $C_i$, $L_i$ and $R_i$, $i=1,\dots,n$, are the capacitance, inductance and resistance of the corresponding capacitors, inductors, and resistors, respectively, and $R_{n+1}$ is the load capacitor. The port-Hamiltonian model of the linear ladder network takes the form
\begin{equation} \label{eq:4.7}
		\dot x = (J_{2n} - R)Q^TQx + u.
\end{equation}
Here $x = (c_1,\phi_1,\dots,c_n,\phi_n)^T$ where $c_i$ and $\phi_i$, for $i=1,\dots,n$, are the charge and the flux of $C_i$ and $L_i$ respectively. $Q$ and $R$ are given as
\begin{equation}
	Q = \text{diag}(C_1^{-1},L_1^{-1},\dots,C_n^{-n},L_n^{-n}), \quad R = \text{diag}(0,R_1,\dots,0,R_n+R_{n+1}),
\end{equation}
$u=(1,0,\dots,0)^T$ is a constant input current and $J_{2n}$ is a skew-symmetric $2n\times 2n$ matrix with -1 and 1 on the superdiagonal and subdiagonal, respectively. 

The energy associated with a port-Hamiltonian system of the form (\ref{eq:4.7}) at time $t$, is given as $H(x(t)) = \frac 1 2 x^T Q^T Q x$. Since $J_{2n}$ is skew symmetric we have that that $\frac d {dt} H(x) = u^T Q^T Q x - x^T Q^T Q R Q^T Q x \leq u^T Q^T Q x$ which is referred to as the \emph{passivity} of the system (\ref{eq:4.7}) \cite{vanderSchaft:1996es,Willems:1972ek}.

Since $J_{2n}$ is full rank, one can always find a coordinate transformation $ x= T \tilde x$ such that $T^{-1} J_{2n} T^{-T} = \mathbb J_{2n}$. The dissipative Hamiltonian formulation of (\ref{eq:4.7}) takes the form
\begin{equation} \label{eq:4.8}
	\dot {\tilde x} = \mathbb J_{2n} \tilde Q^T\tilde Q \tilde x - \tilde Rx + \tilde u,
\end{equation}
where $\tilde Q = QT$, $\tilde R = T^{-1}RT^{-T}Q^TQ$ and $\tilde u = T^{-1} u$. Note that in this case, $\tilde R$ is symmetric and semi-positive definite since $T$ is orthogonal and $R$ is diagonal. The TDD formulation of (\ref{eq:4.8}) takes the form
\begin{equation} \label{eq:4.9}
	\dot{\tilde x} = \mathbb{J}_{2n} \tilde Q^T f(t) + \tilde u, \quad f(t) + \tilde R \int_0^t f(t) = \tilde Q \tilde x.
\end{equation}
The extended Hamiltonian formulation (\ref{eq:2.10.a})-(\ref{eq:2.10.b}) with a quadratic Hamiltonian $H_{\text{ex}}$ can be carried out following Section \ref{sec:2.2}. We note that due to the input $\tilde u$, the Hamiltonian $H_{\text{ex}}$ is time dependent. In fact $\frac{d}{dt} H_{\text{ex}} = \tilde u^T\tilde Q^T \tilde Q \tilde x$. If we define $\overset{\circ}{H}_{\text{ex}} : \mathbb R^{2n}\times \mathcal H^{2n}\times \mathbb R^{2}\to \mathbb R$ as
\begin{equation} \label{eq:4.10}
	\overset{\circ}{H}_{\text{ex}}(\tilde x,\phi,\theta,t,e) = H_{\text{ex}}(\tilde x,\phi,\theta,t) + e, \quad \dot e = - \partial_t H_{\text{ex}},
\end{equation}
it is easily checked that $\frac d {dt} \overset{\circ}{H}_{\text{ex}} =0$ \cite{Hairer:1250576}. However for the time integration of the Hamiltonian system related to $\overset{\circ}{H}_{\text{ex}}$ we can apply a symplectic integrator directly on (\ref{eq:4.10}), since the evolution of $\tilde x$, $\phi$ and $\theta$ does not explicitly depend on $e$. Thus, the passivity of (\ref{eq:4.7}) will be preserved through a symplectic time integration of (\ref{eq:4.9}).

Using an ortho-symplectic reduced basis $A$, the Reduced Dissipative Hamiltonian method can be applied to (\ref{eq:4.9}) to construct a reduced system of the form (\ref{eq:3.12.a})-(\ref{eq:3.12.b}) together with the extended Hamiltonian $\tilde H_{\text{ex}}$. We note that $\frac{d}{dt} \tilde H_{\text{ex}} = (A^+ \tilde u)^T A^T \tilde Q^T \tilde Q A y$, ensuring that the reduced system is passive. Furthermore, the dissipative Hamiltonian structure of the reduced system indicates that the reduced system also carries a port-Hamiltonian structure.

We consider a 100-dimensional ($n=50$) port-Hamiltonian system for the ladder network discussed above. We take $C_i=1$, $L_i = 1$, $R_i=0.2$ for $i=1,\dots,50$, and $R_{51} = 0.4$. We construct the RDH reduced system following Algorithm \ref{alg:3.1}. 

The solution of the RDH method is compared to the main results of \cite{Polyuga:2010gj}, where a passivity-preserving model reduction is developed using a moment matching method at infinity. The charge in $C_1$ is chosen to be the single out put for the moment matching method.

\begin{figure}[t]
\begin{tabular}{cc}
\includegraphics[width=0.5\textwidth]{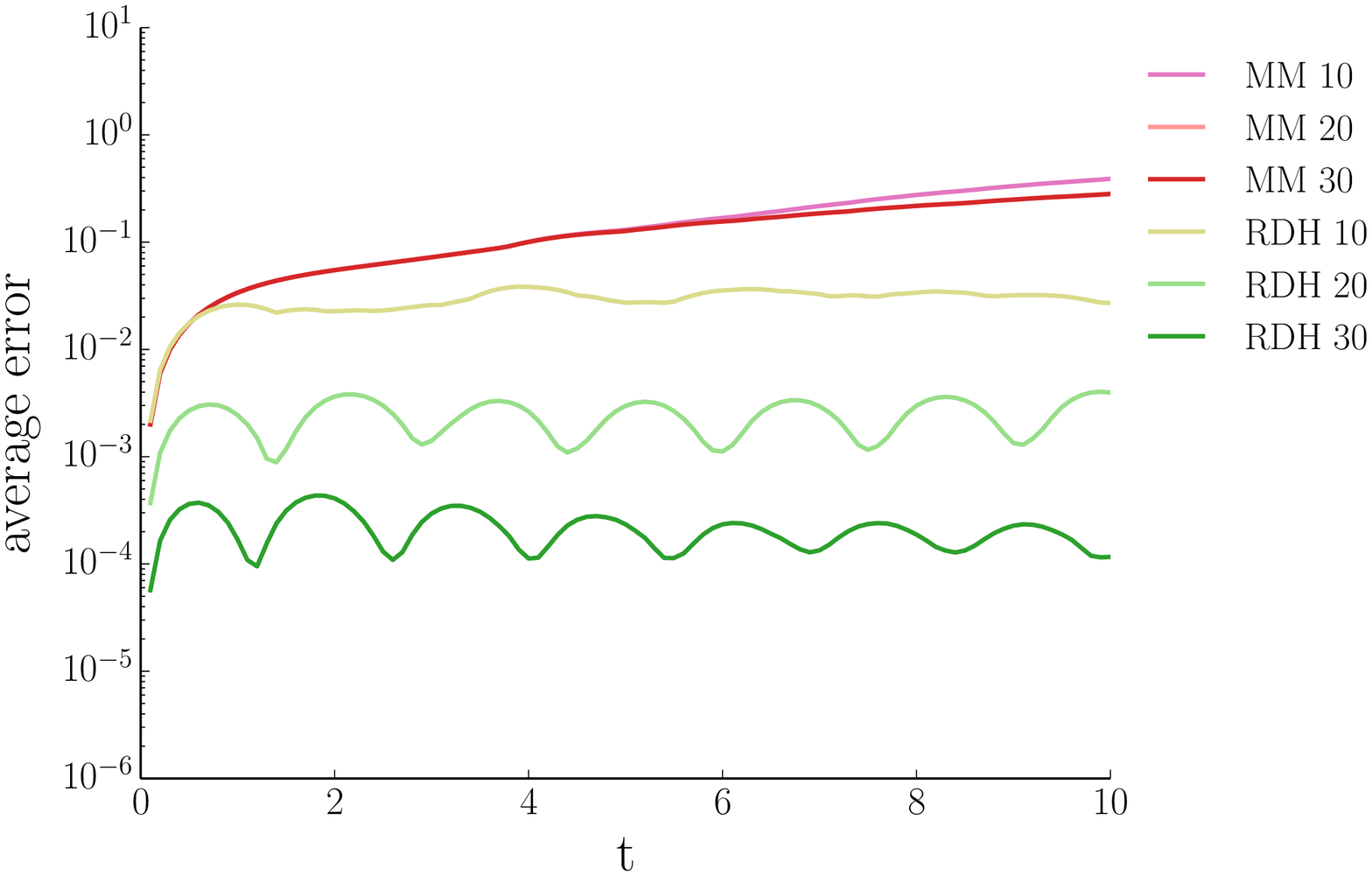} & 
\includegraphics[width=0.5\textwidth]{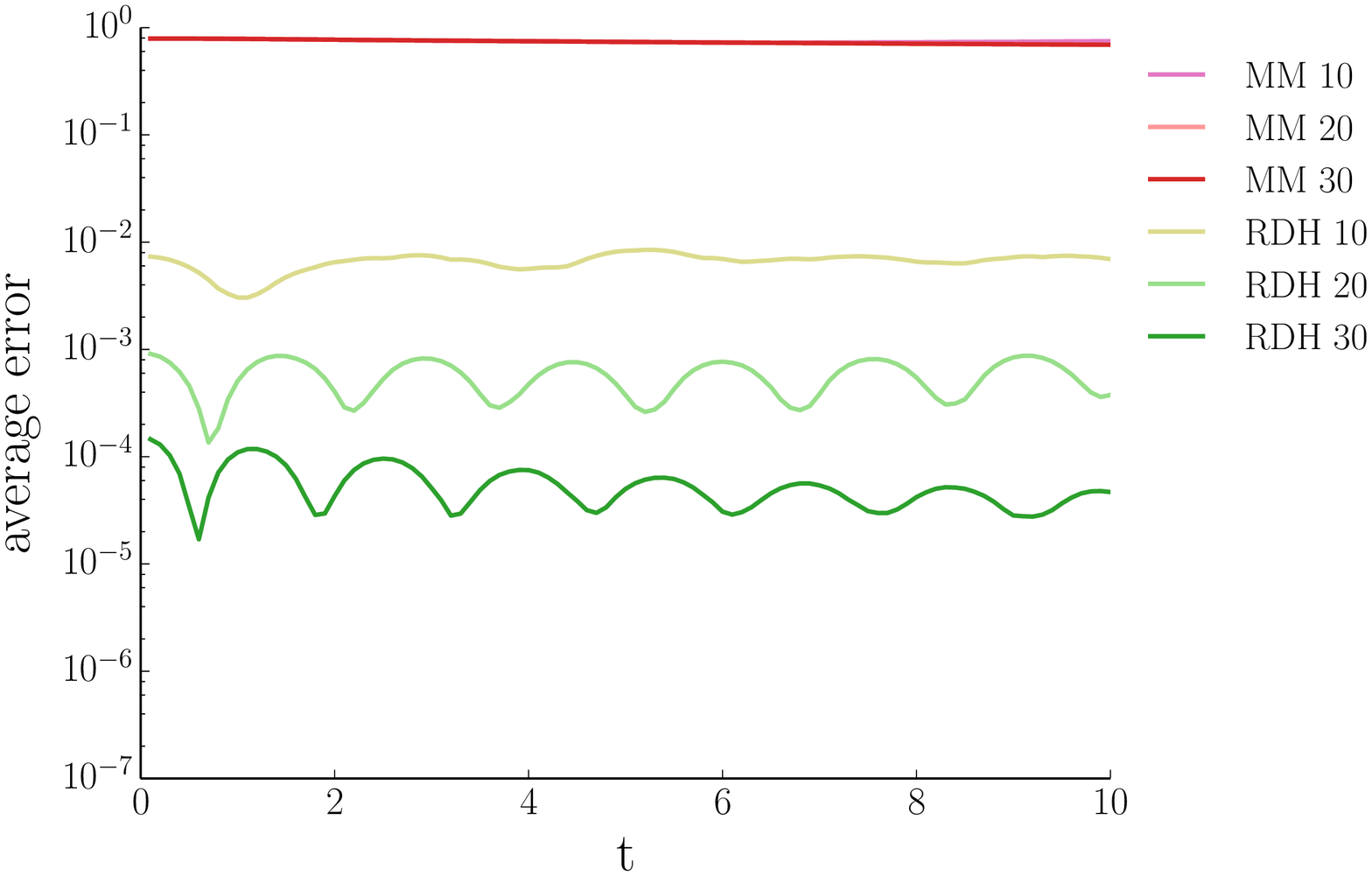} \\
(a) capacitors & (b) inductors \\
\multicolumn{2} {c} {\includegraphics[width=0.5\textwidth]{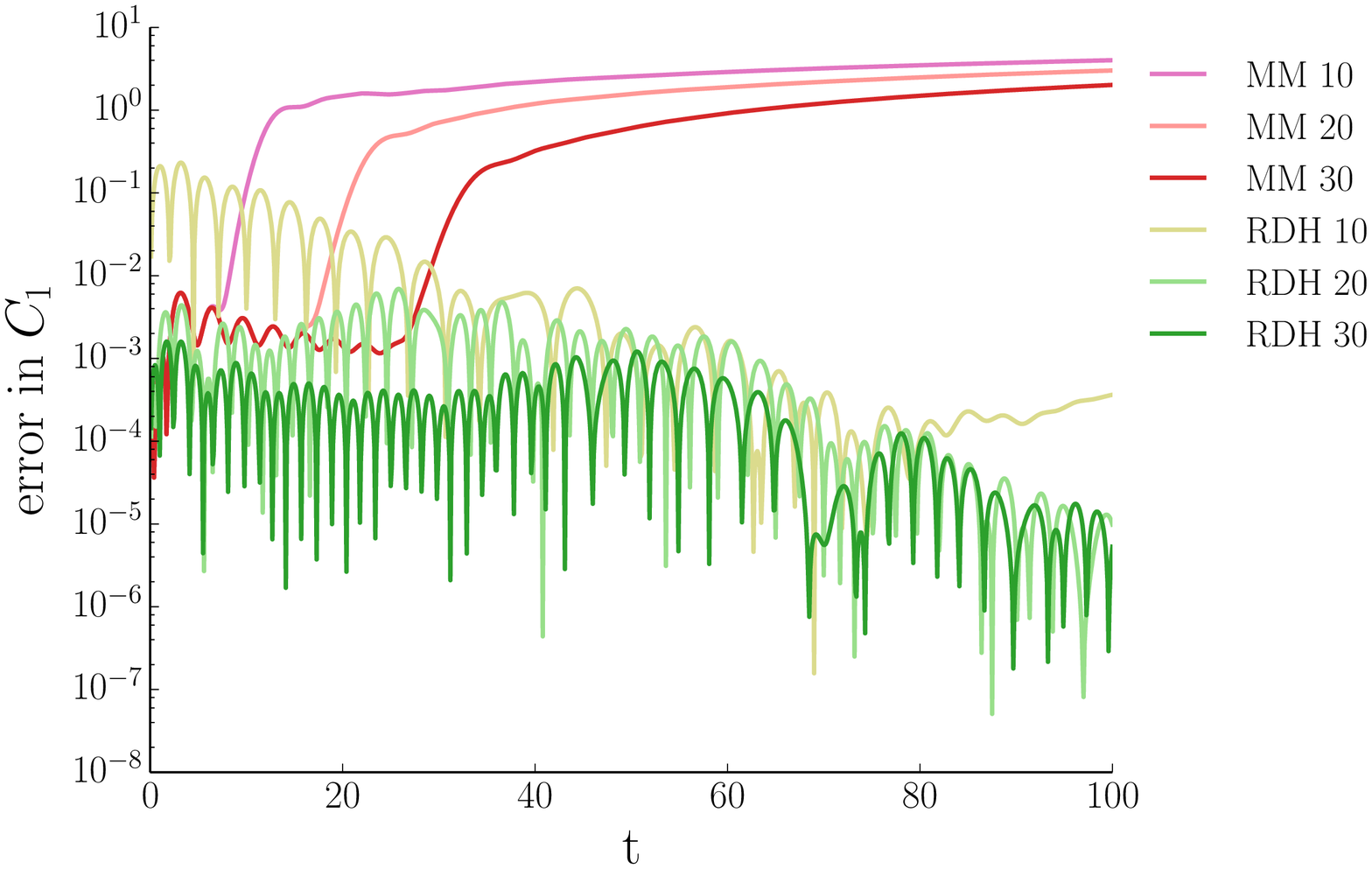}} \\
\multicolumn{2} {c} {(c) charge in $C_1$} 
\end{tabular}
\caption{Error between the full model and the reduced model obtained by the reduced dissipative Hamiltonian method ``RDH'' and the moment matching method ``MM''. (a) The average temporal error of charge in capacitors. (b) the average temporal error of the flux in inductors. (c) The error in $C_1$.} \label{fig:4.3}
\end{figure}

Reduced bases of size $2k = 10$, $2k = 20$ and $2k = 30$ are constructed with the RDH and the moment matching method. Figure \ref{fig:4.3}.(c) shows the error in the charge in $C_1$ for the two methods. We observe that although the moment matching method is bounded over long-time integration, the RDH method provides a significantly more accurate solution. In the moment matching method, the passivity of the reduced system implies that the energy of the system will be bounded by the input energy. However, there is no guarantee that the dissipation of energy in the reduced system mimics the one of the original system. On the other hand, the RDH method allows a correct dissipation of energy through the hidden strings and the symplectic time integration in the RDH method guarantees that the total energy is preserved.

Over short-time integration, we notice that the moment matching method with 10 modes provides a more accurate solution than the RDH with 10 modes. Furthermore, the moment matching method with 20 and 30 modes provide a comparable accuracy to the RDH method with 20 and 30 modes. However, the RDH method maintains the high accuracy during long-time integration, while the moment matching method loses up to 3 orders of magnitude in the accuracy, independent of the number of modes.

Figure \ref{fig:4.3}.(a) and Figure \ref{fig:4.3}.(b) show the average temporal error in the charge and flux of the capacitors and inductors, respectively. The RDH method provides a significantly better accuracy compared to the moment matching method. This is because the charge of $C_1$ is specified as the output of interest in the moment matching method and so it is expected that that method provides low accuracy for computing other outputs. On the other hand, the RDH method not only provides high accuracy in computing the charge for $C_1$ but also high accuracy for all components of the system.

\section{Conclusion} \label{sec:5}

In this paper we present the Reduced Dissipative Hamiltonian method. The method preserves the symplectic structure of dissipative Hamiltonian systems and guarantees the correct dissipation of energy through time integration. The RDH method couples the reduced system with a canonical heat bath such that the reduced system forms a closed system.

The main advantage of the RDH method compared to the existing methods is that it enables the reduced system to be integrated using a symplectic integrator which naturally preserves the Hamiltonian structure and the symplectic symmetry of the Hamiltonian systems. Applying a symplectic integrator to a non-conservative system or using a non-symplectic integrator for the reduced system can cause accumulation of local errors or wrong qualitative solution over long-time integration, respectively.

The numerical simulations illustrate that the RDH method preserves the system energy with significantly higher accuracy than other methods. Furthermore, it is shown that the hidden strings assure that the dissipation of energy in the reduce system mimics the dissipation of energy in the full system. This ensures that the local error do not accumulate over long-time integration.

\bibliographystyle{spmpsci}      
\bibliography{ref}

\end{document}